\documentclass[a4paper, 10pt]{article}

\usepackage{epsfig}
\usepackage{amsmath}
\usepackage{amssymb}
\usepackage{latexsym}
\usepackage{amsfonts}
\usepackage{amsthm}

\usepackage[numbers,sort&compress]{natbib}
\usepackage{graphicx}
\ifx\pdfoutput\undefined  \else
 \fi
\usepackage[centerlast]{caption2}

 \numberwithin{figure}{section}

\begin{document}

\title{
{\Large\bf On the $3$-$\gamma_t$-Critical Graphs of Order
$\Delta(G)+3$} \footnote{The research is supported by Chinese
Natural Science Foundations (60573022). }}
\author{
Haoli Wang,\ Xirong Xu,\ Yang Yuansheng \footnote {corresponding
author's email : yangys@dlut.edu.cn} ,\ Lei Wang   \\\\
Department of Computer Science\\
Dalian University of Technology\\
Dalian, 116024, P. R. China\\  }
\date{}
\maketitle
\begin{abstract}
\noindent Let $\gamma_t(G)$ be the total domination number of graph
$G$, a graph $G$ is $k$-total domination vertex critical (or\ just\
$k$-$\gamma_t$-critical) if $\gamma_t(G)=k$, and for any vertex $v$
of $G$ that is not adjacent to a vertex of degree one,
$\gamma_t(G-v)=k-1$. Mojdeh and Rad \cite{MR06} proposed an open
problem:  Does there exist a $3$-$\gamma_t$-critical graph $G$ of
order $\Delta(G)+3$ with $\Delta(G)$ odd? In this paper, we prove
that there exists a $3$-$\gamma_t$-critical graph $G$ of order
$\Delta(G)+3$ with odd $\Delta(G)\geq 9$.


\end{abstract}

\section{Introduction}
\ \ \ \  Let $G=(V,\ E)$ be a simple connected and undirected graph
with vertex set $V(G)$ and edge set $E(G)$. The {\it open
neighborhood} and the {\it closed neighborhood} of a vertex $v\in
V(G)$ are denoted by $N(v)=\{u \in V(G)\ :\ vu \in E(G)\}$ and
$N[v]=N(v)\cup \{v\}$, respectively. For a vertex set $S \subseteq
V(G)$, $N(S)=\underset{v \in S}\cup N(v)$ and $N[S]=\underset{v \in
S}\cup N[v]$. The {\it degree} of a vertex $v\in V(G)$, denoted by
$deg(v)$, is the number of edges incident with $v$, i.e.
$deg(v)=|N(v)|$.

A set $S \subseteq V(G)$ is a $dominating$ $set$ if for each $v \in
V(G)$ either $v \in S$ or $v$ is adjacent to some $w \in S$. That
is, $S$ is a dominating set if and only if $N[S]=V(G)$. The {\it
domination number} $\gamma(G)$ is the minimum cardinalities of
minimal dominating sets. A set $S \subseteq V(G)$ is a $total$
$dominating$ $set$ if for each $v \in V(G)$, $v$ is adjacent to some
$w \in S$. That is, $S$ is a total dominating set if and only if
$N(S)=V(G)$. The {\it total domination number} $\gamma_t(G)$ is the
minimum cardinalities of minimal total dominating sets. Obviously,
$\gamma(G) \leq \gamma_t(G)$. A $\gamma_t(G)$-set $S$ is a total
dominating set of $G$ with $|S|=\gamma_t(G)$.

A graph $G$ is called {\it vertex domination critical} if
$\gamma(G-v)<\gamma(G)$, for any vertex $v$ of $G$. Further
properties on vertex domination critical graphs were explored in
\cite{AP04,FSW94,FHM95,S90}.

Goddard, Haynes, Henning and van der Merwe \cite{GHHM04} introduced
the concept of total domination vertex critical. A graph $G$ is {\it
$k$-total domination vertex critical} \ (or\ just\
$k$-$\gamma_t$-critical) if $\gamma_t(G)=k,$ and for any vertex $v$
of $G$ that is not adjacent to a vertex of degree one,
$\gamma_t(G-v)=k-1$. For a more detailed treatment of
domination-related parameters and for terminology not defined here,
the reader is referred to \cite{HHS98}.

Mojdeh and Rad \cite{MR06} solved an open problem of
$k$-$\gamma_t$-critical graphs and obtained some results on the
characterization of total domination critical graphs of order
$\Delta(G)+\gamma_t(G)$.

\noindent \textbf{Theorem 1.1.} A cycle $C_n$ is $\gamma_t$-critical
if and only if $n\equiv1,2(\mbox{mod }4)$.

\noindent \textbf{Theorem 1.2.} If $G$ is a $k$-$\gamma_t$-critical
graph of order $\Delta(G)+k$, then $diam(G)=2$.

\noindent \textbf{Theorem 1.3.} There is no $k$-$\gamma_t$-critical
graph of order $\Delta(G)+k$ for $k\geq 4$.

\noindent \textbf{Theorem 1.4.} For any positive integer $k \geq 2$
there exists a $3$-$\gamma_t$-critical graph $G$ of order $2k + 3$
with $\Delta(G)= 2k$.

\noindent \textbf{Theorem 1.5.} There is no $3$-$\gamma_t$-critical
graph $G$ of order $\Delta(G)+3$ with $\Delta(G)= 3,5$.

Where $\Delta(G)\geq 3$ for Theorem 1.2-1.5. They also proposed an
open problem:

Does there exist a $3$-$\gamma_t$-critical graph $G$ of order
$\Delta(G)+3$ with $\Delta(G)$ odd?

In this paper, we solve this open problem by constructing a family
of $3$-$\gamma_t$-critical graphs with  order $\Delta(G)+3$ for odd
$\Delta(G)\geq 9$.

\section{$3$-$\gamma_t$-Critical Graph of order $\Delta(G)+3$}
\indent \indent Let $G$ be a $3$-$\gamma_t$-critical graph of order
$\Delta(G)+3$. Let $x \in V(G)$ be the vertex with
$deg(x)=\Delta(G)$, and $y, z$ be two vertices with $y, z \in
V(G)\setminus N(x)$.

\indent We first give two necessary conditions for
$3$-$\gamma_t$-critical graph of order $\Delta(G)+3$. Let $S_v$ be a
$\gamma_t(G-v)$-set for any $v \in V(G)$. If $S_{v} \cap N(v)\neq
\emptyset$, then $S_{v}$ would be a dominating set of $G$, a
contradiction. Hence we have:

\noindent \textbf{Lemma 2.1.} $S_{v} \bigcap N(v) = \emptyset$ for
any $v \in V(G)$.

\noindent \textbf{Lemma 2.2.} $yz \in E(G)$ and $|N(v)\cap \{ y, z\}
|=1$ for any $v \in N(x)$.
\begin{proof} By Lemma 2.1, $S_{x}=\{y, z\}$, it follows
$yz \in E(G)$, $N(y)\cup N(z)=V(G)\setminus x$, and $|N(v)\cap \{ y,
z\}| \geq 1$ for any $v \in N(x)$. If $|N(v)\cap \{ y, z\}|=2$, then
$S=\{x, v \}$ would be a dominating set of $G$, a contradiction. So
$|N(v)\cap \{ y,z\}|=1$. This completes the proof of Lemma 2.2.
\end{proof}

Then, we proceed to prove that there exists a
$3$-$\gamma_t$-critical graph of order $\Delta(G)+3$ with odd
$\Delta(G) \geq 9$. We formulate such kind of graph $G$ in
Definition 2.1 for the case when $\Delta(G)$ is equivalent to $3$
modulo $4$ and in Definition 2.2 for the case when $\Delta(G)$ is
equivalent to $1$ modulo $4$.

\vspace{11pt}

\noindent \textbf{Definition 2.1.} For $m\geq 3$, let
$G_{4m+2}=(V(G_{4m+2}), E(G_{4m+2}))$, where

\noindent $\begin{array}{llll}
V(G_{4m+2}) =& \{x, y, z\}\cup\\
& \{ y_i \ : 1 \leq i \leq 2m-1\wedge y_i \in N(y)\setminus \{z\} \}\cup\\
& \{ z_i \ : 1 \leq i \leq 2m\wedge z_i \in N(z)\setminus \{y\}\}\\
\end{array}$

\noindent \ \ and

\noindent $\begin{array}{llll}
E(G_{4m+2}) =& \{ xy_i \ : 1 \leq i \leq 2m-1 \}\cup\\
& \{ y_iy \ : 1 \leq i \leq 2m-1 \}\cup\\
& \{ xz_i \ : 1 \leq i \leq 2m \}\cup\\
& \{ z_iz \ : 1 \leq i \leq 2m \}\cup\\
& \{ y_iy_{1+(i+m-2) \mbox { mod }(2m-1)} \ : 1 \leq i \leq 2m-1 \}\cup\\
& \{ z_iz_j \ : 1 \leq i < j \leq 2m \wedge j \neq i+m \}\cup\\
& \{ y_iz_j \ : 1 \leq i \leq 2m-1 \wedge 1 \leq j \leq 2m \wedge j \neq i \wedge j \neq i+1 \}.\\
\end{array}$

Then, $\Delta(G_{4m+2}) = 4m-1$, and the order of $G_{4m+2}$ is
$\Delta(G_{4m+2})+3$.

As an example, $G_{14}$ and its complement graph $\overline {
G_{14}}$ are shown in Figure 2.1.

\begin{figure}[h]
\centering
\includegraphics{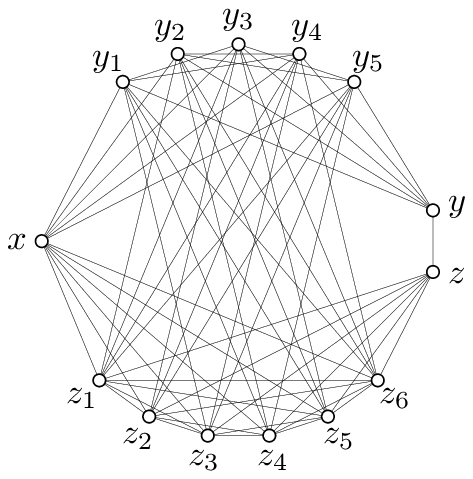}
\includegraphics{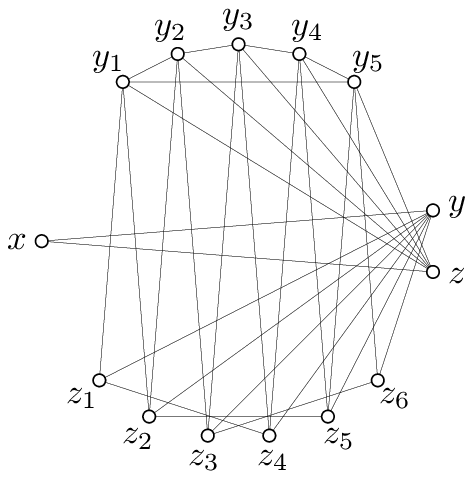}
\caption{$G_{14}$ and its complement graph $\overline {G_{14}}$.}
\end{figure}

\noindent \textbf{Lemma 2.3.} $\gamma_t(G_{4m+2})=3$.

\begin{proof}(a) Since $S = \{ x,y_1,y \}$ is a total
dominating set of $G_{4m+2}$ with $|S| = 3$, we have
$\gamma_t(G_{4m+2}) \leq 3$.

(b) We prove $\gamma_t(G_{4m+2}) \geq 3$ by showing that $S=\{u,v\}$
does not dominate $G_{4m+2}$ for any edge $uv\in E(G_{4m+2})$ in the
following cases:

\noindent Case 1. $S\cap (N[y]\setminus \{z\})=\emptyset$. Then $y$
is not dominated by $S$.

\noindent Case 2. $S\cap (N[z]\setminus \{y\})=\emptyset$. Then $z$
is not dominated by $S$.

\noindent Case 3. $S=\{y, z\}$. Then $x$ is not dominated by $S$.

\noindent Case 4. $S=\{y_i, z_j\}$, $1 \leq i \leq 2m-1$, $1\leq
j\leq 2m$. By symmetry, we only need to consider $1 \leq j \leq m$.
Since $y_j$ is not dominated by $z_j$, then $y_iy_j\in E(G_{4m+2})$.

\noindent Case 4.1. $1 \leq j \leq m-1$. Then $i \in \{j+m-1, j+m
\}$, $z_{j+m}$ is not dominated by $S$.

\noindent Case 4.2. $j=m$. Then $i \in \{1, 2m-1\}$. If $i=1$, then
$y_{m-1}$ is not dominated by $S$. If $i=2m-1$, then $z_{2m}$ is not
dominated by $S$.

 By (a) and (b), we have $\gamma_t(G_{4m+2})=3$.
\end{proof}

\noindent \textbf{Lemma 2.4.} $\gamma_t(G_{4m+2}-v)=2$ for any
vertex $v\in V(G_{4m+2})$.

\begin{proof} We prove this lemma by showing that there exists a subset $S_v\subseteq V(G)$ with $|S_v|=2$ for
any $v\in V(G_{4m+2})$.

\noindent Case 1. $v=x$. Let $S_v=\{y,z\}$.

\noindent Case 2. $v=y$. Let $S_v=\{x, z_1 \}$.

\noindent Case 3. $v=z$. Let $S_v=\{x, y_1 \}$.

\noindent Case 4. $v=y_i$, $1 \leq i \leq 2m-1$. If $i \leq m$, let
$S_v=\{y_{i+m-2},z_{i} \}$. If $i \geq m+1$, let
$S_v=\{y_{i-m+2},z_{i+1} \}$.

\noindent Case 5. $v=z_i$, $1 \leq i \leq 2m$. By symmetry, we only
need to consider $1 \leq i \leq m$. Let $S_v=\{y_i, z_{i+m}\}$.

From the above cases, we have that $\gamma_t(G_{4m+2}-v)=2$ for any
vertex $v\in V(G_{4m+2})$. \end{proof}

As an immediate result of Lemma 2.3 and Lemma 2.4, we have the
following:

\noindent \textbf{Theorem 2.5.} $G_{4m+2}$ is a
$3$-$\gamma_t$-critical graph.

\vspace{11pt}

\noindent \textbf{Definition 2.2.} For $m\geq 3$, let
$G_{4m}=(V(G_{4m}), E(G_{4m}))$, where

\noindent $\begin{array}{llll}
V(G_{4m}) =& \{x, y, z\}\cup\\
& \{ y_i \ : 1 \leq i \leq 2m-2\wedge y_i \in N(y)\setminus \{z\} \}\cup\\
& \{ z_i \ : 1 \leq i \leq 2m-1\wedge z_i \in N(z)\setminus
y\}\\
\end{array}$

\noindent \ \ and

\noindent $\begin{array}{llll}
E(G_{4m}) =& \{ xy_i \ : 1 \leq i \leq 2m-2 \}\cup\\
& \{ y_iy \ : 1 \leq i \leq 2m-2 \}\cup\\
& \{ xz_i \ : 1 \leq i \leq 2m-1 \}\cup\\
& \{ z_iz \ : 1 \leq i \leq 2m-1 \}\cup\\
& \{ y_iy_j \ : 1 \leq i < j \leq 2m-2 \wedge j \neq i+m-1 \}\cup\\
& \{ z_iz_{1+(i+m-2) \mbox{ mod } (2m-1)} \ : 1 \leq i \leq 2m-1 \}\cup\\
& \{ y_iz_j \ : 1 \leq i \leq 2m-2 \wedge 1 \leq j \leq 2m-1 \wedge j \neq i \wedge j \neq i+1 \}\cup\\
& \{ y_{1}z_{2} \}.\\
\end{array}$

Then, $\Delta(G_{4m})=4m-3$, and the order of $G_{4m}$ is
$\Delta(G_{4m})+3$.

As an example, $G_{12}$ and its complement graph $\overline {
G_{12}}$ are shown in Figure 2.2.

\begin{figure}[h]
\centering
\includegraphics{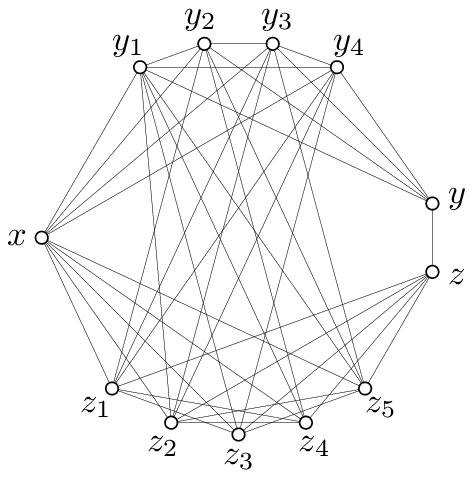}
\includegraphics{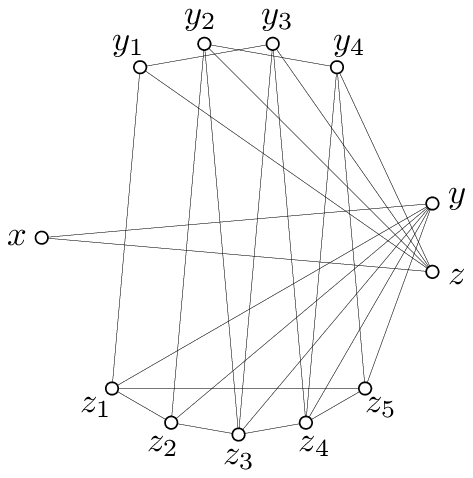}
\caption{$G_{12}$ and its complement graph $\overline {G_{12}}$.}
\end{figure}

\noindent \textbf{Lemma 2.6.} $\gamma_t(G_{4m})=3$.

\begin{proof}(a) Since $S = \{ x,y_1,y \}$ is a total
dominating set of $G_{4m}$ with $|S| = 3$, we have $\gamma_t(G_{4m})
\leq 3$.

(b) We prove $\gamma_t(G_{4m}) \geq 3$ by showing that $S=\{u,v\}$
does not dominate $G_{4m}$ for any edge $uv\in E(G_{4m})$ in the
following cases:

\noindent Case 1. $S\cap (N[y]\setminus \{z\})=\emptyset$. Then $y$
is not dominated by $S$.

\noindent Case 2. $S\cap (N[z]\setminus \{y\})=\emptyset$. Then $z$
is not dominated by $S$.

\noindent Case 3. $S=\{y, z\}$. Then $x$ is not dominated by $S$.

\noindent Case 4. $S=\{y_i, z_j\}$,$1 \leq i \leq 2m-2$, $1 \leq j
\leq 2m-1$.

\noindent Case 4.1. $1 \leq i \leq m-1$. Since $z_i$ is not
dominated by $y_i$, it follows $z_iz_j \in E(G_{4m})$, $j \in
\{i+m-1, i+m\}$, $y_{i+m-1}$ is not dominated by $S$.

\noindent Case 4.2. $m \leq i \leq 2m-2$. Since $z_{i+1}$ is not
dominated by $y_i$, it follows $z_{i+1}z_j \in E(G_{4m})$, $j \in
\{i-m+1, i-m+2\}$. If $m+1 \leq i \leq 2m-2$ or $i=m$ and $j=1$,
then $y_{i-m+1}$ is not dominated by $S$. If $i=m$ and $j=2$, then
$z_m$ is not dominated by $S$.

By (a) and (b), we have $\gamma_t(G_{4m})=3$. \end{proof}

\noindent \textbf{Lemma 2.7.} $\gamma_t(G_{4m}-v)=2$ for any vertex
$v\in V(G_{4m})$.

\begin{proof} We prove this lemma by showing that there exists a subset $S_v\subseteq V(G)$ with $|S_v|=2$ for
any $v\in V(G_{4m})$.

\noindent Case 1. $v=x$. Let $S_v=\{y,z\}$.

\noindent Case 2. $v=y$. Let $S_v=\{x, z_1\}$.

\noindent Case 3. $v=z$. Let $S_v=\{x, y_1\}$.

\noindent Case 4. $v=y_i$, $1 \leq i \leq 2m-2$. If  $1 \leq i \leq
m-1$, let $S_v=\{y_{i+m-1}, z_i\}$. If  $m \leq i \leq 2m-2$, let
$S_v=\{y_{i-m+1}, z_{i+1}\}$.

\noindent Case 5. $v=z_{i}$, $1 \leq i \leq 2m-1$.

If $1 \leq i \leq m$, let $S_v=\{y_{i},z_{1+(i+m) \mbox{ mod }
(2m-1)} \}$. If $m+1 \leq i \leq 2m-1$, let
$S_v=\{y_{i-1},z_{1+(i+m-3) \mbox{ mod }(2m-1)} \}$.

From the above cases, we have that $\gamma_t(G_{4m}-v)=2$ for any
vertex $v\in V(G_{4m})$.
\end{proof}

As an immediate result of Lemma 2.6 and Lemma 2.7, we have the
following:

\noindent \textbf{Theorem 2.8.} $G_{4m}$ is a
$3$-$\gamma_t$-critical graph.

\vspace{10pt}

From Theorem 2.5 and Theorem 2.8, we conclude that there exists a
$3$-$\gamma_t$-critical graph $G$ of order $\Delta(G)+3$ for an
arbitrary odd integer $\Delta(G)\geq 9$. By theorem 1.4, we know
that there is no $3$-$\gamma_t$-critical graph $G$ of order
$\Delta(G)+3$ with $\Delta(G)= 3,5$. Thus, in the rest of this
paper, we shall prove that there does not exist a
$3$-$\gamma_t$-critical graph of order $\Delta(G)+3$ with
$\Delta(G)=7$.

Denote $ N(y)\setminus \{z\}=\{y_1,y_2,\cdots,y_{|N(y)\setminus
\{z\}|}\}$ and $ N(z)\setminus
\{y\}=\{z_1,z_2$,\\$\cdots,z_{|N(z)\setminus \{y\}|}\}.$

\noindent \textbf{Observation 2.9.} \\
\indent (a) there exists a pair of vertices $(y_i,z_j) (1\leq i\leq
|N(y)\setminus \{z\}|,1\leq j\leq |N(z)\setminus \{y\}|)$, such that
$y_iz_j\in
E(G)$; \\
\indent (b) there exists at least one vertex $z_j (1\leq j\leq
|N(z)\setminus \{y\}|)$, such that $y_iz_j\not \in E(G)$ for each
vertex $y_i (1\leq i\leq
|N(y)\setminus \{z\}|)$; \\
\indent (c) there exists at least one vertex $y_i (1\leq i\leq
|N(y)\setminus \{z\}|)$, such that $y_iz_j\not \in E(G)$ for each
vertex $z_j (1\leq j\leq |N(z)\setminus \{y\}|)$.

\noindent \textbf{Theorem 2.10.} There is no $3$-$\gamma_t$-critical
graph of order $\Delta(G)+3$ with $\Delta(G)=7$.

\begin{proof} Without loss of generality, by Lemma 2.2, we may assume
that $|N(y)|<|N(z)|$ and distinguish the possible
$3$-$\gamma_t$-critical graphs of order $\Delta(G)+3$ with
$\Delta(G)=7$ into three cases:

\noindent Case 1. $|N(y)\setminus \{z\}|=1$ and $|N(z)\setminus
\{y\}|=6$.

By Observation 2.9, we derive a contradiction.

\vspace{10pt}

\noindent Case 2. $|N(y)\setminus \{z\}|=2$ and $|N(z)\setminus
\{y\}|=5$.

Without loss of generality, we may assume that $|N(y_1)\cap
N(z)|\geq |N(y_2)\cap N(z)|$.

\noindent Case 2.1. $|N(y_1)\cap(N(z)\setminus \{y\})|=4$, say
$y_1z_1, y_1z_2, y_1z_3, y_1z_4 \in E(G)$. Then $S_{z_1}=\{ y_2,
z_5\}$, $y_1y_2, z_2z_5, z_3z_5, z_4z_5 \in E(G)$, and $\{ y_1,
z_2\}$ would be a dominating set of $G$(see Figure 2.3. $G_{10.1}$),
a contradiction.

\noindent Case 2.2. $|N(y_1)\cap(N(z)\setminus \{y\})|=3$, say
$y_1z_1, y_1z_2, y_1z_3\in E(G)$. Then $S_{z_1} \in \{ \{ y_2,
z_4\}, \{y_2, z_5\} \}$, say $S_{z_1}=\{ y_2, z_4\}$, $y_1y_2,
z_2z_4, z_3z_4 \in E(G)$, $S_{z_2}=\{ y_2, z_5\}$, $z_1z_5, z_3z_5
\in E(G)$, and $\{ y_1, z_3\}$ would be a dominating set of $G$(see
Figure 2.3. $G_{10.2}$), a contradiction.

\begin{figure}[h]
\renewcommand{\captionlabeldelim}{:}
\centering
\includegraphics[scale=0.58]{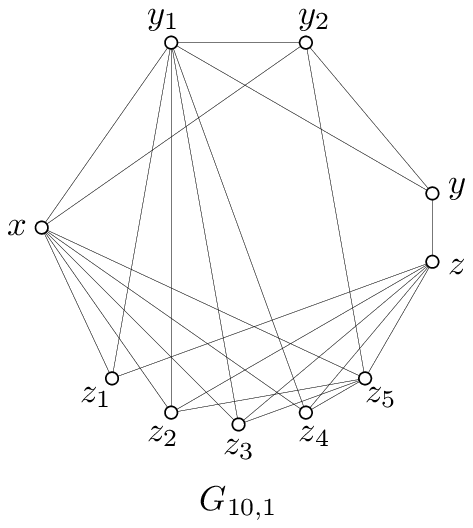}
\includegraphics[scale=0.58]{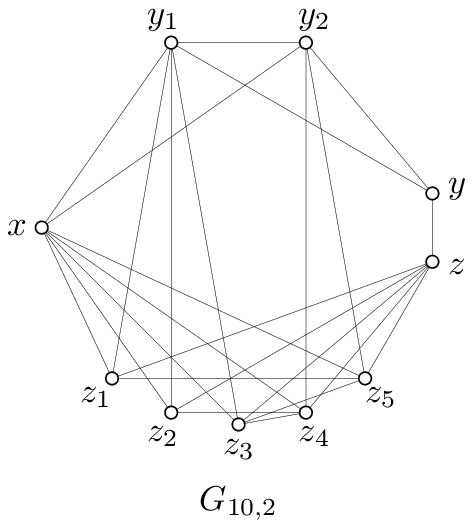}
\includegraphics[scale=0.58]{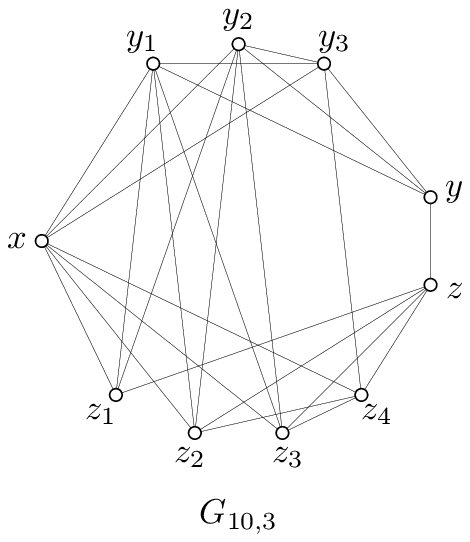}
\includegraphics[scale=0.58]{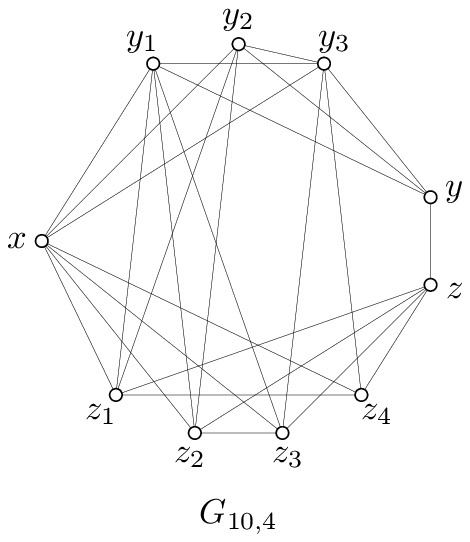}
\caption{$G_{10.1}-G_{10.4}$ for non-$3$-$\gamma_t$-critical graph
$G$ of order $\Delta(G)+3$ with $\Delta(G)=7$.}
\end{figure}

\noindent Case 3. $|N(y)\setminus \{z\}|=3$ and $|N(z)\setminus
\{y\}|=4$.

Without loss of generality, we may assume that $|N(y_1)\cap
N(z)|\geq |N(y_2)\cap N(z)|$ and $|N(y_1)\cap N(z)|\geq |N(y_3)\cap
N(z)|$.

\vspace{10pt}

\noindent Case 3.1. $|N(y_1)\cap (N(z)\setminus \{y\})|=3$, say
$y_1z_1,y_1z_2,y_1z_3\in E(G)$. We may assume that $|N(y_2)\cap
\{z_1,z_2,z_3\}|\geq |N(y_3)\cap \{z_1,z_2,z_3\}|$.

\noindent Case 3.1.1. $|N(y_2)\cap \{z_1,z_2,z_3\}|=3$. Then
$z_1,z_2,z_3\not \in N(y_3)$, $S_{z_1}=\{y_3,z_4\}$. It follows that
$y_1y_3,y_2y_3,z_2z_4,z_3z_4\in E(G)$, $\{y_1,z_2\}$ would be a
dominating set of $G$(see Figure 2.3. $G_{10.3}$), a contradiction.

\noindent Case 3.1.2. $|N(y_2)\cap \{z_1,z_2,z_3\}|=2$, say
$y_2z_1,y_2z_2\in E(G)$. Then $z_1,z_2\not \in N(y_3)$,
$S_{z_1}\in\{\{y_3,z_3\},\{y_3,z_4\}\}$. If $S_{z_1}=\{y_3,z_3\}$,
then $y_2y_3,z_2z_3\in E(G)$, $S_{z_2}=\{y_3,z_4\}$,
$y_1y_3,z_1z_4\in E(G)$, and $\{y_1,z_1\}$ would be a dominating set
of $G$(see Figure 2.3. $G_{10.4}$), a contradiction. If
$S_{z_1}=\{y_3,z_4\}$, then $y_1y_3,z_2z_4\in E(G)$, and
$\{y_1,z_2\}$ would be a dominating set of $G$(see Figure 2.4.
$G_{10.5}$), a contradiction.

\noindent Case 3.1.3. $|N(y_2)\cap \{z_1,z_2,z_3\}|=1$, say
$y_2z_1\in E(G)$.

\noindent Case 3.1.3.1. $|N(y_3)\cap \{z_1,z_2,z_3\}|=1$, say
$y_3z_2\in E(G)$. Then $S_{z_1}\in\{\{y_3,z_2\},\{y_3,z_4\}\}$.

\noindent Case 3.1.3.1.1. $S_{z_1}=\{y_3,z_2\}$, then
$y_2y_3,z_2z_3\in E(G)$, $z_4$ is adjacent to at least one vertex of
$\{y_3,z_2\}$, $S_{z_2}\in\{\{y_2,z_1\},\{y_2,z_4\}\}$.

\noindent Case 3.1.3.1.1.1. $S_{z_2}=\{y_2,z_1\}$. Then $z_1z_3\in
E(G)$, $S_{z_3}\in \{\{y_2,z_4\},$ $\{y_3,z_4\}\}$. If
$S_{z_3}=\{y_2,z_4\}$, then $y_1y_2,z_2z_4\in E(G)$, and
$\{y_1,z_2\}$ would be a dominating set of $G$(see Figure 2.4.
$G_{10.6}$), a contradiction. If $S_{z_3}=\{y_3,z_4\}$, then
$y_1y_3,z_1z_4\in E(G)$, and $\{y_1,z_1\}$ would be a dominating set
of $G$(see Figure 2.4. $G_{10.7}$), a contradiction.

\noindent Case 3.1.3.1.1.2. $S_{z_2}=\{y_2,z_4\}$. Then
$y_1y_2,z_3z_4\in E(G)$, $S_{z_3}=\{y_2,z_1\}$, $z_1z_2\in E(G)$,
and $\{y_3,z_2\}$ would be a dominating set of $G$(see Figure 2.4.
$G_{10.8}$), a contradiction.
\begin{figure}[h]
\renewcommand{\captionlabeldelim}{:}
\centering
\includegraphics[scale=0.58]{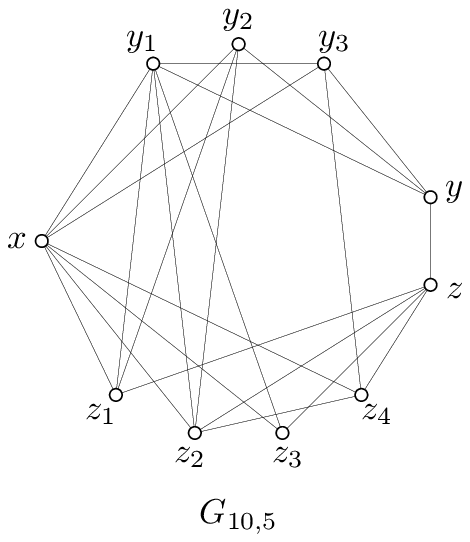}
\includegraphics[scale=0.58]{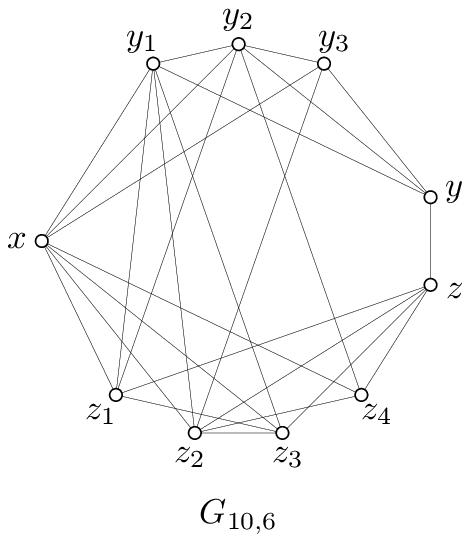}
\includegraphics[scale=0.58]{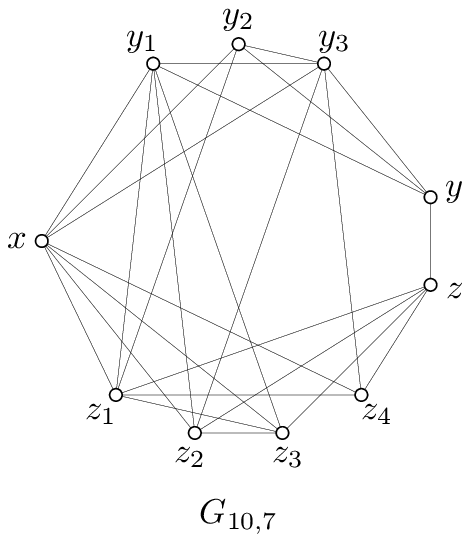}
\includegraphics[scale=0.58]{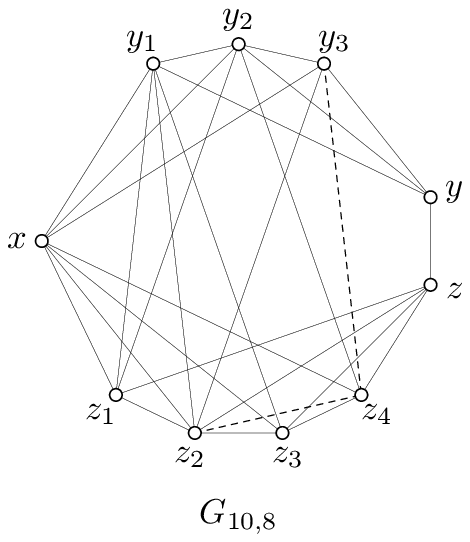}
\includegraphics[scale=0.58]{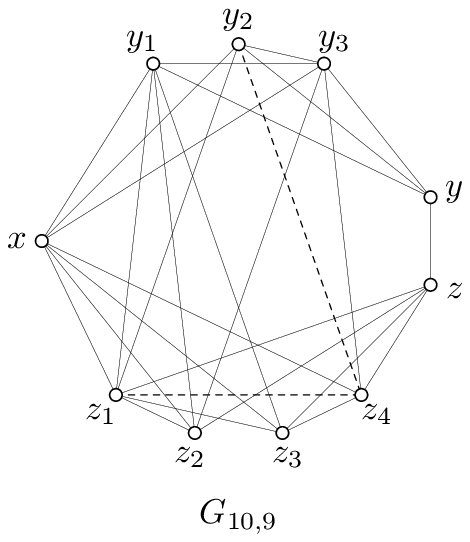}
\includegraphics[scale=0.58]{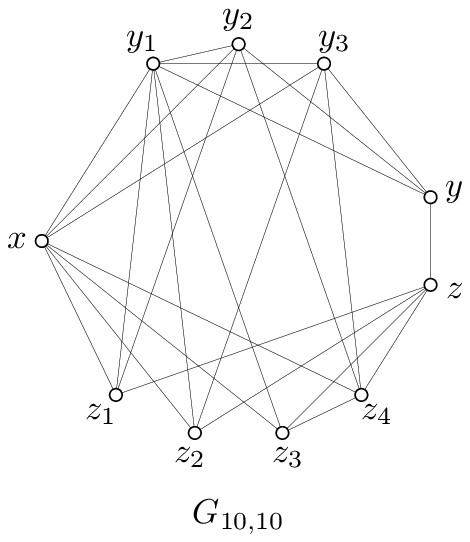}
\includegraphics[scale=0.58]{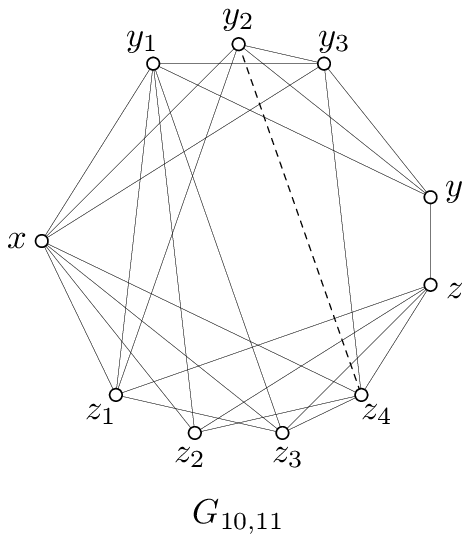}
\includegraphics[scale=0.58]{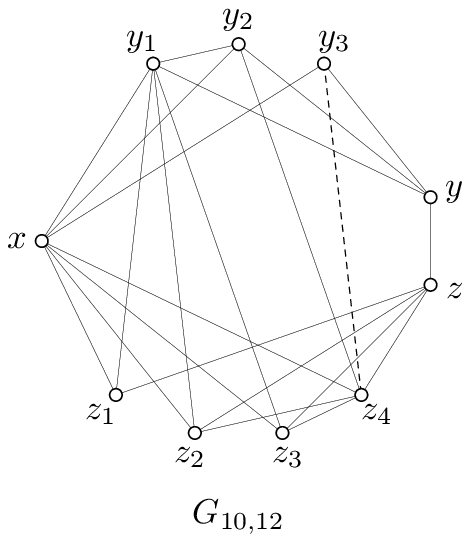}
\caption{$G_{10.5}-G_{10.12}$ for non-$3$-$\gamma_t$-critical graph
$G$ of order $\Delta(G)+3$ with $\Delta(G)=7$.}
\end{figure}

\noindent Case 3.1.3.1.2. $S_{z_1}=\{y_3,z_4\}$, then
$y_1y_3,z_3z_4\in E(G)$, $S_{z_2}\in\{\{y_2,z_1\},$ $\{y_2,z_4\}\}$.

\noindent Case 3.1.3.1.2.1. $S_{z_2}=\{y_2,z_1\}$. Then
$y_2y_3,z_1z_3\in E(G)$, $z_4$ is adjacent to at least one vertex of
$\{y_2,z_1\}$, $S_{z_3}=\{y_3,z_2\}$, $z_1z_2\in E(G)$, and
$\{y_2,z_1\}$ would be a dominating set of $G$(see Figure 2.4.
$G_{10.9}$), a contradiction.

\noindent Case 3.1.3.1.2.2. $S_{z_2}=\{y_2,z_4\}$. Then $y_1y_2\in
E(G)$, and $\{y_1,z_3\}$ would be a dominating set of $G$(see Figure
2.4. $G_{10.10}$), a contradiction.

\noindent Case 3.1.3.2. $|N(y_3)\cap \{z_1,z_2,z_3\}|=0$. Then
$S_{z_1}=\{y_3,z_4\}$, $y_1y_3,z_2z_4,$ $z_3z_4\in E(G)$,
$S_{z_2}=\{y_2,z_1\}$, $y_2y_3,z_1z_3\in E(G)$, and $S_{z_3}\in
\{\{y_1,z_1\},$ $\{y_1,z_2\},\{y_2,z_1\},\{y_2,z_4\},\{y_3,z_4\}\}$
would be a dominating set of $G$(see Figure 2.4. $G_{10.11}$), a
contradiction.

\noindent Case 3.1.4. $|N(y_2)\cap \{z_1,z_2,z_3\}|=0$. Then
$S_{z_1}\in \{\{y_2,z_4\},\{y_3,z_4\}\}$, say $S_{z_1}=\{y_2,z_4\}$,
then $y_1y_2,z_2z_4,z_3z_4\in E(G)$, $S_{z_2}\in
\{\{y_1,z_1\},\{y_1,z_3\},$ $\{y_2,z_4\}, \{y_3,z_4\}\}$ would be a
dominating set of $G$(see Figure 2.4. $G_{10.12}$), a contradiction.

\vspace{11pt}

\noindent Case 3.2. $|N(y_1)\cap (N(z)\setminus \{y\})|=2$, say
$y_1z_1,y_1z_2\in E(G)$.  We may assume that $|N(y_2)\cap
\{z_1,z_2\}|\geq |N(y_3)\cap \{z_1,z_2\}|$.

\noindent Case 3.2.1. $|N(y_2)\cap \{z_1,z_2\}|=2$. Then
$z_1,z_2\not \in N(y_3)$, $S_{y_1}\in \{\{y_3,z_3\},$
$\{y_3,z_4\}\}$, say $S_{y_1}=\{y_3,z_3\}$. It follows that
$y_2y_3,z_1z_3,z_2z_3\in E(G)$, $S_{z_1}=\{y_3,z_4\}$,
$y_1y_3,z_2z_4\in E(G)$, and $\{y_1,z_2\}$ would be a dominating set
of $G$(see Figure 2.5. $G_{10.13}$), a contradiction.

\noindent Case 3.2.2. $|N(y_2)\cap \{z_1,z_2\}|=1$, say $y_2z_1\in
E(G)$. Then $z_1\not \in N(y_3)$.

\noindent Case 3.2.2.1. $y_3z_2\in E(G)$. We may assume that
$|N(y_2)\cap \{z_3,z_4\}|\geq |N(y_3)\cap \{z_3,z_4\}|$. There are
four subcases.

\noindent Case 3.2.2.1.1. $y_2z_3,y_3z_3\in E(G)$. Then $S_{z_1}\in
\{\{y_3,z_2\},\{y_3,z_3\}\}$, say $S_{z_1}=\{y_3,z_2\}$. It follows
$y_2y_3,z_2z_4\in E(G)$, $S_{z_2}\in \{\{y_2,z_1\},\{y_2,z_3\}\}$.

\noindent Case 3.2.2.1.1.1. $S_{z_2}=\{y_2,z_1\}$. It follows
$z_1z_4\in E(G)$. $S_{z_4}\in\{\{y_2,z_3\},$ $\{y_3,z_3\} \}$, say
$\{y_2,z_3\}$, then $y_1y_2,z_2z_3\in E(G)$, and $\{y_1,z_2\}$ would
be a dominating set of $G$(see Figure 2.5. $G_{10.14}$), a
contradiction.
\begin{figure}[h]
\renewcommand{\captionlabeldelim}{:}
\centering
\includegraphics[scale=0.58]{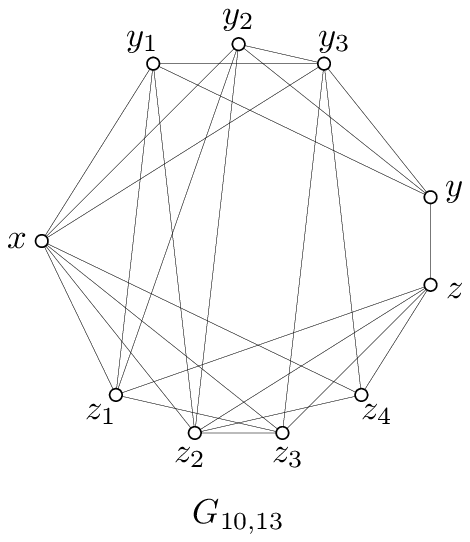}
\includegraphics[scale=0.58]{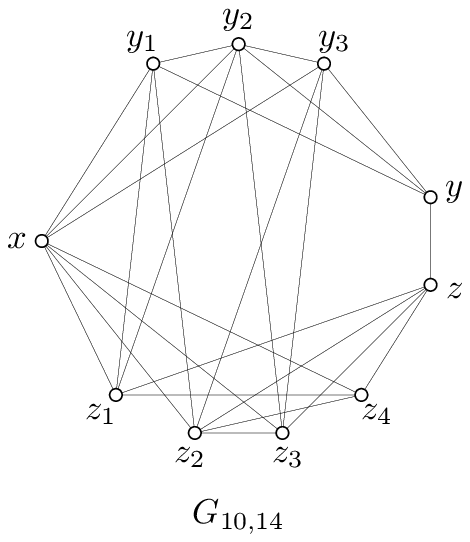}
\includegraphics[scale=0.58]{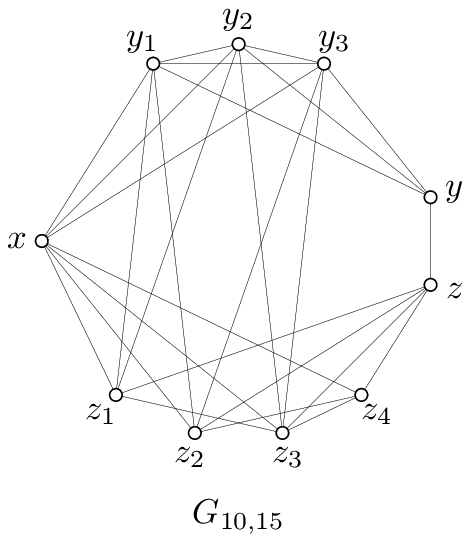}
\includegraphics[scale=0.59]{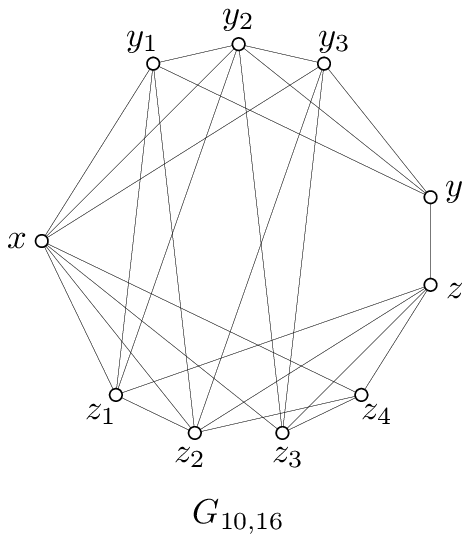}
\caption{$G_{10.13}-G_{10.16}$ for non-$3$-$\gamma_t$-critical graph
$G$ of order $\Delta(G)+3$ with $\Delta(G)=7$.}
\end{figure}

\noindent Case 3.2.2.1.1.2. $S_{z_2}=\{y_2,z_3\}$. It follows
$y_1y_2,z_3z_4\in E(G)$, $S_{z_4}\in\{\{y_1,z_1\}$, $\{y_2,z_1\}\}$.
If $S_{z_4}=\{y_1,z_1\}$, then $y_1y_3,z_1z_3\in E(G)$, and
$\{y_3,z_3\}$ would be a dominating set of $G$(see Figure 2.5.
$G_{10.15}$), a contradiction. If $S_{z_4}=\{y_2,z_1\}$, then
$z_1z_2\in E(G)$, and $\{y_3,z_2\}$ would be a dominating set of
$G$(see Figure 2.5. $G_{10.16}$), a contradiction.

\noindent Case 3.2.2.1.2. $y_2z_3,y_3z_4\in E(G)$. Then $S_{y_1}\in
\{\{y_2,z_3\},\{y_3,z_4\}\}$, say $S_{y_1}=\{y_2,z_3\}$. It follows
$y_2y_3,z_2z_3,z_3z_4\in E(G)$, $S_{z_3}=\{y_1,z_1\}$, $y_1y_3,$
$z_1z_4\in E(G)$, and $\{y_3,z_4\}$ would be a dominating set of
$G$(see Figure 2.6. $G_{10.17}$), a contradiction.

\noindent Case 3.2.2.1.3. $y_2z_3\in E(G)$ and $y_3$ is not adjacent
to any vertex of $\{z_3,z_4\}$. Then $S_{z_1}=\{y_3,z_2\}$,
$y_2y_3,z_2z_3,z_2z_4\in E(G)$. It follows $S_{z_2}=\{y_2,z_1\}$,
$z_1z_4\in E(G)$, $S_{z_4}=\{y_2,z_3\}$, $y_1y_2\in E(G)$, and
$\{y_1,z_2\}$ would be a dominating set of $G$(see Figure 2.6.
$G_{10.18}$), a contradiction.

\noindent Case 3.2.2.1.4. Both $y_2$ and $y_3$ are not adjacent to
any vertex of $\{z_3,z_4\}$. Then $S_{y_1}\in
\{\{y_2,z_1\},\{y_3,z_2\}\}$ would be a dominating set of $G$(see
Figure 2.6. $G_{10.19}$), a contradiction.

\noindent Case 3.2.2.2. $y_3z_2\not \in E(G)$. Then $S_{z_1}\in
\{\{y_3,z_3\},\{y_3,z_4\}\}$, say $S_{z_1}=\{y_3,z_3\}$. It follows
$y_1y_3,z_2z_3\in E(G)$, $y_2$ is adjacent to at least one vertex of
$\{y_3,z_3\}$ and $z_4$ is adjacent to at least one vertex of
$\{y_3,z_3\}$.

\noindent Case 3.2.2.2.1. $y_2z_3,z_3z_4\in E(G)$. Then
$S_{z_3}=\{y_1,z_1\}$. It follows $z_1z_4\in E(G)$,
$S_{z_4}=\{y_1,z_2\}$, $y_1y_2\in E(G)$, and $\{y_2,z_3\}$ would be
a dominating set of $G$(see Figure 2.6. $G_{10.20}$), a
contradiction.
\begin{figure}[hr]
\renewcommand{\captionlabeldelim}{:}
\centering
\includegraphics[scale=0.59]{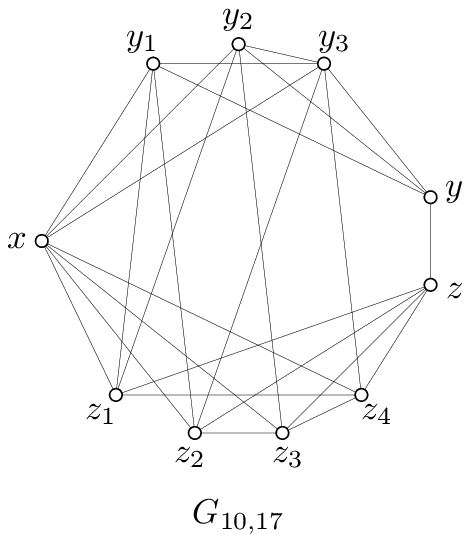}
\includegraphics[scale=0.59]{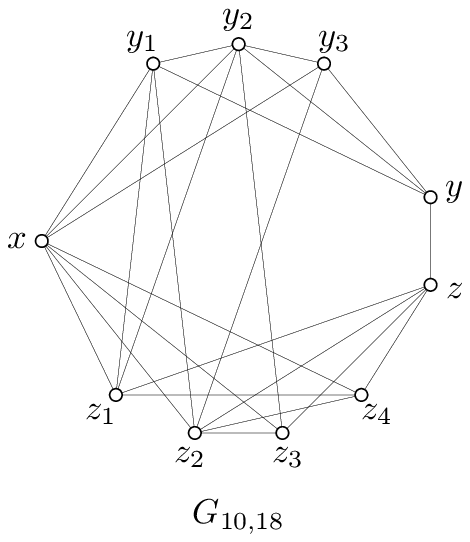}
\includegraphics[scale=0.59]{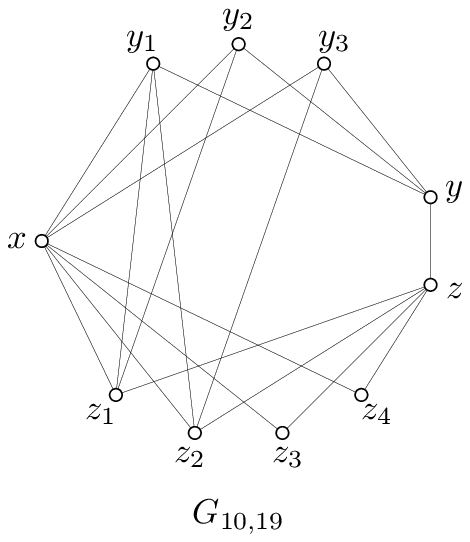}
\includegraphics[scale=0.59]{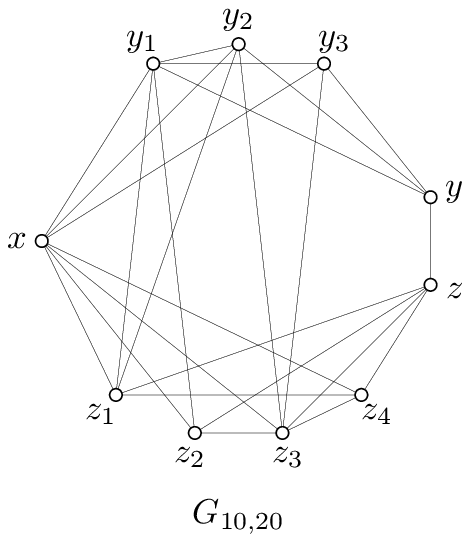}
\includegraphics[scale=0.59]{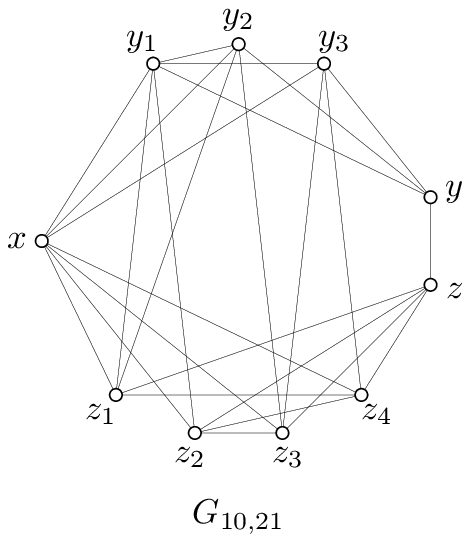}
\includegraphics[scale=0.59]{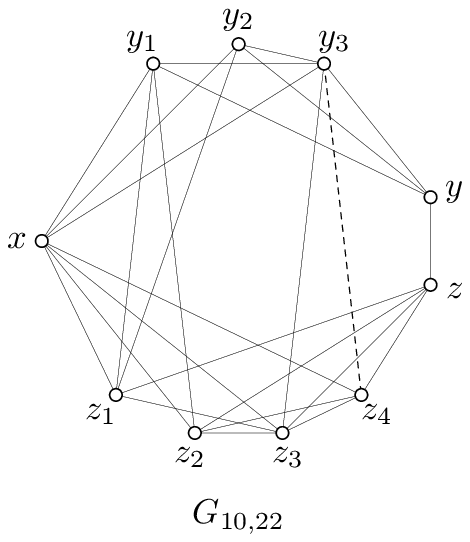}
\includegraphics[scale=0.59]{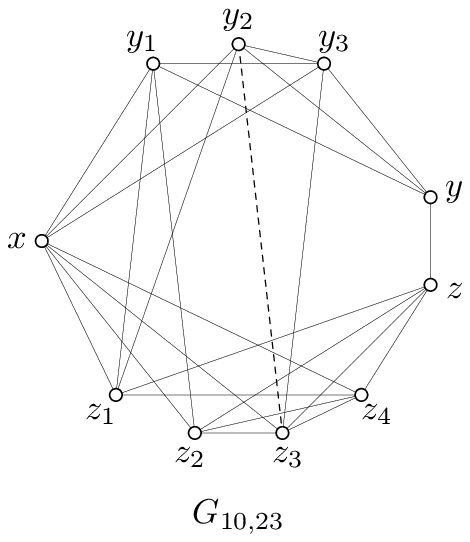}
\includegraphics[scale=0.59]{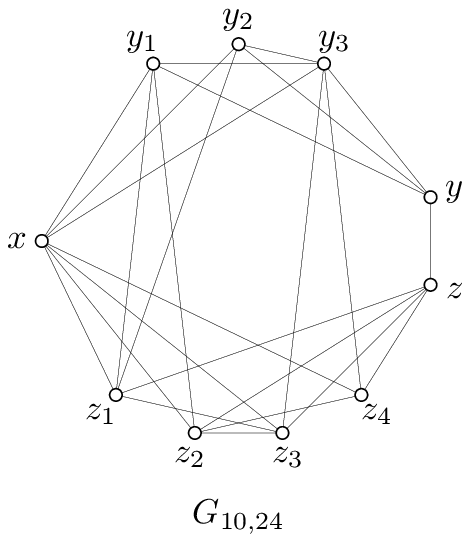}
\caption{$G_{10.17}-G_{10.24}$ for non-$3$-$\gamma_t$-critical graph
$G$ of order $\Delta(G)+3$ with $\Delta(G)=7$.}
\end{figure}

\noindent Case 3.2.2.2.2. $y_2z_3,y_3z_4\in E(G)$. Then
$S_{z_3}=\{y_1,z_1\}$. It follows $z_1z_4\in E(G)$, $S_{z_4}\in
\{\{y_1,z_2\},\{y_2,z_3\}\}$, then $y_1y_2\in E(G)$,
$S_{y_2}=\{y_3,z_4\}$, $z_2z_4\in E(G)$, and $\{y_1,z_2\}$ would be
a dominating set of $G$(see Figure 2.6. $G_{10.21}$), a
contradiction.

\noindent Case 3.2.2.2.3. $y_2y_3,z_3z_4\in E(G)$. Then
$S_{y_2}=\{y_1,z_2\}$. It follows $z_2z_4\in E(G)$,
$S_{z_2}=\{y_2,z_1\}$, $z_3$ is adjacent to at least one vertex of
$\{y_2,z_1\}$ and $z_4$ is adjacent to at least one vertex of
$\{y_2,z_1\}$. Since $|N(y_1)\cap N(z)|\geq |N(y_2)\cap N(z)|$,
$y_2$ is adjacent to at most one vertex of $\{z_3,z_4\}$. Hence
$z_1$ is adjacent to at least one vertex of $\{z_3,z_4\}$. If
$z_1z_3\in E(G)$, then $S_{z_3}\in
\{\{y_1,z_1\},\{y_1,z_2\},\{y_2,z_1\},\{y_2,z_4\},\{y_3,z_4\}\}$
would be a dominating set of $G$(see Figure 2.6. $G_{10.22}$), a
contradiction. If $z_1z_4\in E(G)$, then $S_{z_4}\in
\{\{y_1,z_1\},\{y_1,z_2\},\{y_2,z_1\}$, $\{y_2,z_3\},\{y_3,z_3\}\}$
would be a dominating set of $G$(see Figure 2.6. $G_{10.23}$), a
contradiction.

\noindent Case 3.2.2.2.4. $y_2y_3,y_3z_4\in E(G)$. Then
$S_{y_2}=\{y_1,z_2\}$. It follows $z_2z_4\in E(G)$,
$S_{z_2}=\{y_2,z_1\}$, $z_1$ is adjacent to at least one vertex of
$\{z_3,z_4\}$. If $z_1z_3\in E(G)$, then $\{y_3,z_3\}$ would be a
dominating set of $G$(see Figure 2.6. $G_{10.24}$), a contradiction.
If $z_1z_4\in E(G)$, then $\{y_3,z_4\}$ would be a dominating set of
$G$(see Figure 2.7. $G_{10.25}$), a contradiction.

\noindent Case 3.2.3. $|N(y_2)\cap \{z_1,z_2\}|=0$. Then $S_{z_1}\in
\{\{y_2,z_3\},\{y_2,z_4\},\{y_3,z_3\},$ $\{y_3,z_4\}\}$, say
$S_{z_1}=\{y_2,z_3\}$. It follows that $y_1y_2,z_2z_3\in E(G)$,
$S_{y_2}=\{y_3,z_4\}$, $y_1y_3,z_1z_4,z_2z_4\in E(G)$, and
$\{y_1,z_2\}$ would be a dominating set of $G$(see Figure 2.7.
$G_{10.26}$), a contradiction.

\begin{figure}[hr]
\renewcommand{\captionlabeldelim}{:}
\centering
\includegraphics[scale=0.59]{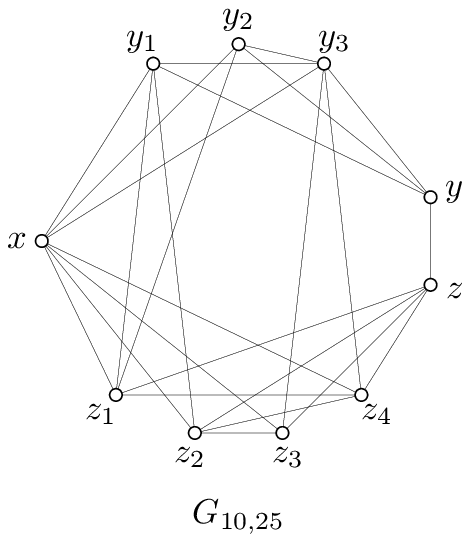}
\includegraphics[scale=0.59]{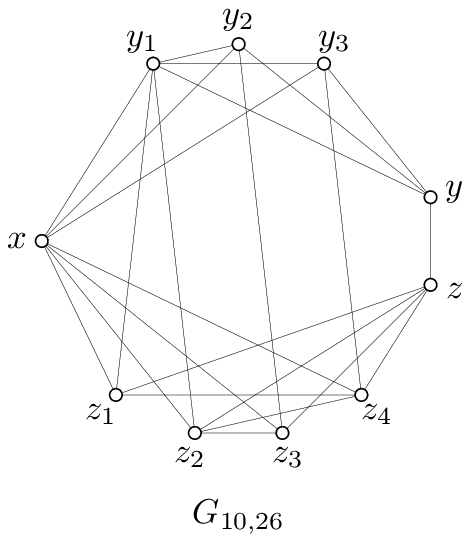}
\includegraphics[scale=0.59]{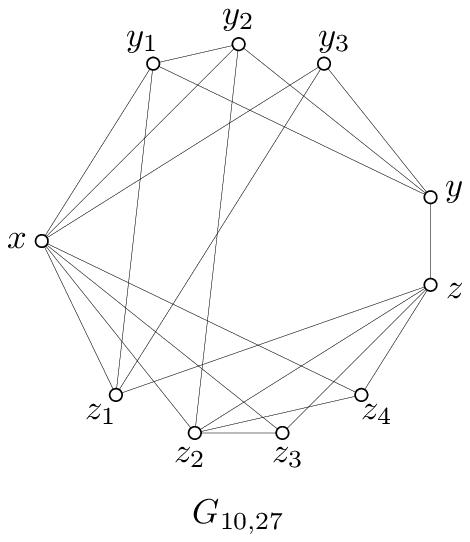}
\includegraphics[scale=0.59]{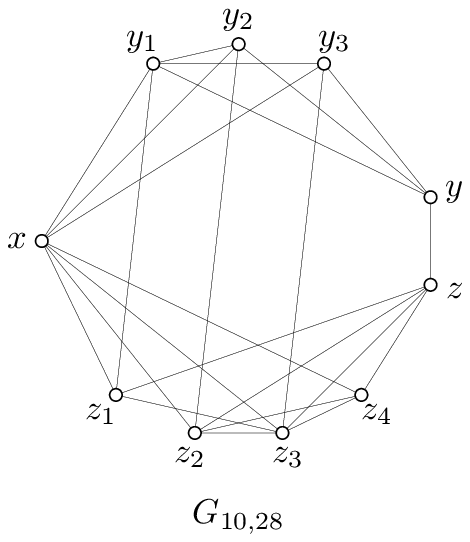}
\caption{$G_{10.25}-G_{10.28}$ for non-$3$-$\gamma_t$-critical graph
$G$ of order $\Delta(G)+3$ with $\Delta(G)=7$.}
\end{figure}

\noindent Case 3.3. $|N(y_1)\cap (N(z)\setminus \{y\})|=1$, say
$y_1z_1\in E(G)$. Then $S_{z_1}\in
\{\{y_2,z_2\},\{y_2,z_3\},\{y_2,z_4\},\{y_3,z_2\},\{y_3,z_3\},\{y_3,z_4\}\}$,
say $S_{z_1}=\{y_2,z_2\}$. It follows $y_1y_2,z_2z_3,z_2z_4\in
E(G)$. If $y_3z_1\in E(G)$, $S_{y_1}\in \{\{y_2, z_2 \},\{y_3, z_1
\} \}$ would be a dominating set of $G$(see Figure 2.7.
$G_{10.27}$), a contradiction. If $y_3$ is adjacent to one vertex of
$\{z_3,z_4\}$, say $z_3$, then $S_{y_2} =\{y_3,z_3\}$,
$y_1y_3,z_1z_3$, $z_3z_4\in E(G)$, and $S_{z_3}\in
\{\{y_1,z_1\},\{y_2,z_2\}\}$ would be a dominating set of $G$(see
Figure 2.7. $G_{10.28}$), a contradiction.

This completes the proof of Theorem 2.10.
\end{proof}

Using Theorem 1.1, 1.3, 1.4, 2.5, 2.8 and 2.10, we come to the
conclusion of this paper:

\noindent \textbf{Theorem 2.11.} For any $\Delta(G) \geq 2$, there
exists a $3$-$\gamma_t$-critical graph $G$ of order $\Delta(G)+3$ if
and only if $\Delta(G) \neq 3, 5, 7$.


\begin{thebibliography}{99}

\bibitem{AP04}
N. Ananchuen and M. D. Plummer,
\newblock Matching properties in domination critical graphs,
\newblock {\it Discrete Math.}, 277 (2004) 1-13.

\bibitem{FSW94}
O. Favaron, D. Sumner and E. Wojcicka,
\newblock The diameter of domination critical graphs,
\newblock {\it J.Graph Theory}, 18 (1994) 723-734.

\bibitem{FHM95}
J. Fulman, D. Hanson, G. MacGillivray,
\newblock Vertex domination-critical graphs,
\newblock {\it Networks}, 25 (1995) 41-43.

\bibitem{GHHM04}
W. Goddard, T. W. Haynes, M. A. Henning and L. C. van der Merwe,
\newblock The diameter of total domination vertex critical graphs,
\newblock {\it Discrete Math.}, 286 (2004) 255-261.

\bibitem{HHS98}
T. W. Haynes, S. T. Hedetniemi, P. J. Slater,
\newblock Fundamentals of domination in graphs,
\newblock {\it Marcel Dekker, Inc., NewYork}, 1998.

\bibitem{MR06}
D. A. Mojdeh, N. J. Rad,
\newblock On the total domination critical graphs,
\newblock {\it Electronic Notes in Discrete Math.}, 24 (2006) 89-92.

\bibitem{S90}
D. P. Sumner,
\newblock Critical concepts in domination,
\newblock {\it Discrete Math.}, 86 (1990) 33-46.

\end{thebibliography}
\end{document}